\documentclass[11pt]{amsart}
\usepackage[T1]{fontenc}
\usepackage{lmodern}
\usepackage{microtype}
\usepackage{amsmath,amssymb,amsthm}
\usepackage{mathtools}
\mathtoolsset{showonlyrefs}

\usepackage[foot]{amsaddr}
\usepackage[hang,flushmargin]{footmisc}

\usepackage[hyphens]{url}
\usepackage[hidelinks]{hyperref}
\usepackage{doi}

\usepackage[a4paper]{geometry}
\newtheorem{theorem}{Theorem}[section]
\newtheorem{corollary}[theorem]{Corollary}
\newtheorem{lemma}[theorem]{Lemma}
\newtheorem{proposition}[theorem]{Proposition}
\newtheorem{remark}[theorem]{Remark}
\newtheorem{definition}[theorem]{Definition}
\newcommand{\R}{\mathbb R}
\newcommand{\Z}{\mathbb Z}
\newcommand{\N}{\mathbb N}

\newcommand{\sepdist}{\operatorname{sep}}

\newcommand{\vol}{\operatorname{vol}}

\title{Uniqueness of the hexagonal lattice's universal optimality}

\author[M.~Faulhuber]{Markus Faulhuber}
\email{markus.faulhuber@univie.ac.at}
\address{Faculty of Mathematics, University of Vienna \\ Oskar-Morgenstern-Platz 1, 1090 Vienna, Austria}

\author[P.~Petersen]{Philipp Petersen}
\email{philipp.petersen@univie.ac.at}
\address{Faculty of Mathematics, University of Vienna \\ Kolingasse 14-16, 1090 Vienna, Austria}

\thanks{This research was supported in whole or in part by the Austrian Science Fund (FWF) \href{https://doi.org/10.55776/PAT5102224}{10.55776/PAT5102224}. For open access purposes, the authors have applied a CC BY public copyright license to any author accepted manuscript version arising from this submission.}



\date{}

\usepackage{xcolor}

\begin{document}

\begin{abstract}
    We build on a recent proof of the universal optimality of the planar hexagonal lattice. We extend this result to demonstrate that, up to Euclidean motion, the hexagonal lattice is the only periodic configuration with this property.
\end{abstract}

\maketitle

\vspace{-\baselineskip}
\section{Energy minimization and universal optimality}
We study the problem of arranging points in periodic configurations to minimize certain energies. Towards this goal, we call a function $p: \R_{>0} \to \R$ a potential. For a discrete set $X \subset \R^n$, $n\in \N$, we denote the number of points in $X$ by $|X|$. For a point configuration $C \subset \R^n$, we define its $p$-energy~as
\begin{equation}
    E_p(C) = \lim_{R \to \infty} \frac{1}{|C \cap \mathbf{B}_R|} \sum_{\substack{\gamma, \gamma' \in C \cap \mathbf{B}_R \\ \gamma \neq \gamma'}} p(\lVert \gamma-\gamma' \rVert),
\end{equation}
where $\mathbf{B}_R = \{ x \in \R^n \mid \lVert x \rVert < R\}$ denotes the open ball of radius $R>0$ centered at the origin. The limit may not exist, in which case we define the energy to be infinite. The quantity
\begin{equation}
    \rho(C) = \lim_{R \to \infty} \frac{|C \cap \mathbf{B}_R|}{\vol(\mathbf{B}_R)}
\end{equation}
is called the density of the point configuration $C$.

The energy minimization problem is as follows. For a fixed potential and a fixed density, we seek to find the minimizer of $E_p$ among all point configurations. First, we formulate this problem for Gaussian energies $\varphi_\alpha: r \mapsto e^{-\alpha r^2}$, $\alpha > 0$. Later, in Subsection \ref{sec:generalPotentials} we describe how the problem can naturally be extended to so-called completely monotone functions of squared distance.
\begin{definition}
    A point configuration $C_0$ with $\rho(C_0) = 1$ is called universally optimal if we have
    \begin{equation}
        E_{\varphi_\alpha}(C_0) \leq E_{\varphi_\alpha}(C), \quad \text{ for all }  C \text{ with } \rho(C)=1 \text{ and all } \alpha > 0.
    \end{equation}
\end{definition}
If $C_0$ is universally optimal, then for every $\rho>0$, the rescaled configuration $\widetilde{C}_0 = \rho^{-1/d}C_0$ minimizes every Gaussian energy among configurations of density $\rho$, i.e., 
\begin{equation}\label{eq:universal-optimality-general-rho}
    E_{\varphi_\alpha}(\widetilde{C}_0) \leq E_{\varphi_\alpha}(C), \quad \text{ for all } C \text{ with } \rho(C)=\rho \text{ and } \alpha > 0.
\end{equation}

In dimension two, a celebrated result of Montgomery \cite{Montgomery1988} establishes that the hexagonal (or triangular) lattice minimizes every Gaussian energy among \emph{lattices} of fixed density. The hexagonal lattice of density one is defined as
\begin{equation}\label{eq:hexagonal-lattice}
    \Lambda_2 = \sqrt{\tfrac{2}{\sqrt{3}}}
    \begin{pmatrix}
        1 & \frac{1}{2}\\
        0 & \frac{\sqrt{3}}{2}
    \end{pmatrix}
    \Z^2.
\end{equation}
Cohn and Kumar conjectured that the hexagonal lattice is universally optimal among periodic point configurations of fixed density in $\R^2$ \cite[Conjecture~9.4]{CohnKumar2007}. This universal optimality was recently proved by OpenAI \cite{OpenAI}. 

A natural associated question is whether $\Lambda_2$ is, up to Euclidean motion, the only point configuration that has this property. To make this question precise, we note that uniqueness among \emph{arbitrary} configurations does not hold. 
This is because deleting finitely many points from $\Lambda_2$ preserves its density and Gaussian energies, but produces a configuration that is not congruent to $\Lambda_2$. We therefore restrict our uniqueness question to periodic point configurations.

It is explicitly mentioned in \cite{OpenAI} that only minimality of $\Lambda_2$ for all Gaussian energies is shown, but not the uniqueness of $\Lambda_2$ for the considered problem.

Our contribution is to close this gap by proving the following result.
\begin{theorem}\label{thm:main-uniqueness}
    For every periodic point configuration $C \subset \R^2$ with $\rho(C) = 1$ we have
    \begin{equation}\label{eq:energyminimality}
        E_{\varphi_\alpha}(\Lambda_2) \leq E_{\varphi_\alpha} (C), \quad \text{ for all } \alpha > 0.
    \end{equation}
    Moreover, the following statements apply.
    \begin{enumerate}
        \item[(A)] If equality in \eqref{eq:energyminimality} holds for every $\alpha>0$, then $C$ is the hexagonal lattice $\Lambda_2$ up to Euclidean motion.
        \item[(B)] If equality holds in \eqref{eq:energyminimality} for at least one $\alpha>0$, then $C$ is the hexagonal lattice $\Lambda_2$ up to Euclidean motion.
    \end{enumerate}
\end{theorem}

Statement \eqref{eq:energyminimality} was shown in \cite{OpenAI}. Clearly, statement (B) implies (A). We nevertheless give a separate, simpler proof of (A), for which we only need the asymptotic decay of Gaussian energy and uniqueness of the densest periodic sphere packing in the plane. No further results from \cite{OpenAI} are used. Establishing (B) additionally uses properties of the auxiliary functions constructed in \cite{OpenAI}, which we identify below.

\subsection{General potentials}\label{sec:generalPotentials}
Using standard techniques, many of the arguments around universal optimality can naturally be extended from Gaussian potentials to completely monotone functions of squared distance.
\begin{definition}
    A function $f:\R_{>0} \to \R$ is \emph{completely monotone} if $f\in C^\infty(\R_{>0})$ and
    \begin{equation}
        (-1)^{k} f^{(k)}(s) \geq 0, \quad \text{ for all } s \in \R_{>0}, \ k=0,1,2,\ldots.
    \end{equation}
    A potential $g:\R_{>0}\to\R$ is \emph{completely monotone of squared distance} if $g(r)=f(r^2)$ for some completely monotone function $f$.
\end{definition}

Observe that completely monotone functions are in particular nonnegative, nonincreasing, and convex. 
The Bernstein--Widder theorem characterizes completely monotone functions via the Laplace transform of nonnegative Borel measures \cite{Bernstein1929,Widder1941}.
\begin{theorem}[Bernstein-Widder] \label{thm:bernstein-widder}
    Let $f: \R_{>0} \to \R_{\geq 0}$ be completely monotone. Then there exists a non-negative Borel measure $\mu_f$ such that
    \begin{equation}
        f(s) = \int_0^\infty e^{-\alpha s} \, d\mu_f(\alpha).
    \end{equation}
\end{theorem}
Consequently, every completely monotone potential of squared distance has the representation
\begin{equation}
    g(r)=f(r^2)=\int_0^\infty e^{-\alpha r^2}\,d\mu_f(\alpha).
\end{equation}
Thus, completely monotone functions of squared distance are nonnegative, continuous superpositions of Gaussians $\varphi_\alpha$ and the constant function. Therefore, inequality \eqref{eq:energyminimality} remains valid with $\varphi_\alpha$ replaced by any completely monotone potential of squared distance. 

Moreover, the following uniqueness result for the hexagonal lattice, similar to Theorem \ref{thm:main-uniqueness}, holds in this case. 
\begin{corollary}\label{cor:cm}
    Let $g\not\equiv0$ be a completely monotone potential of squared distance, and suppose $E_g(\Lambda_2)<\infty$. Then
    \begin{equation}
        E_g(C)\ge E_g(\Lambda_2)
    \end{equation}
    for every periodic configuration $C\subset\R^2$, $\rho(C) = 1$, with equality if and only if $C$ is congruent to~$\Lambda_2$.
\end{corollary}

\subsection{Related work}

Universal optimality was shown to hold for the integers $\Z$ in dimension 1 \cite{CohnKumar2007} and, in seminal work by Cohn, Kumar, Miller, Radchenko, and Viazovska \cite{CKMRV2022}, in dimensions 8 and 24, where the minimizers are the $\mathsf{E}_8$-lattice and the Leech lattice $\Lambda_{24}$, respectively. In all cases, the minimizer is also, up to Euclidean motion, the unique periodic configuration of density 1 that yields the densest sphere packing in that dimension. While sphere packing in one dimension is trivial, the seminal breakthrough by Viazovska \cite{Viazovska2017} solving the sphere packing problem in dimension 8, soon followed by a proof for dimension 24 \cite{CKMRV2017}, is based on the linear programming bounds established by Cohn and Elkies \cite{CohnElkies2003}. In all proved cases, including OpenAI's result in dimension 2 \cite{OpenAI}, universal optimality could be established by constructing a magic function saturating the linear programming bounds for energy minimization \cite{CohnKumar2007}. Yet, it should be noted that an alternative, earlier proof in dimension 1 was given by Ventevogel and Nijboer \cite{VentevogelNijboer1979}.

The problem of identifying optimal configurations for energy minimization in the plane has a rich history deeply intertwined with the Epstein zeta function. For a two-dimensional lattice $\Lambda \subset \mathbb{R}^2$, the function is defined for $\text{Re}(s) > 1$ by 
\begin{equation}
    \zeta_\Lambda(s) = \sum_{\lambda \in \Lambda \setminus \{0\}} \frac{1}{\lVert \lambda \rVert^{2s}},
\end{equation}
and extended via analytic continuation to $\mathbb{C} \setminus \{1\}$. For real $s>1$, this sum equals the normalized energy $E_{p_s}(\Lambda)$ for the Riesz potential $p_s(r)=r^{-2s}$. For $0<s\le1$, the positive interaction sum diverges, so the analytically continued zeta values must be distinguished from this energy.

Rankin~\cite{Rankin1953} proved that the hexagonal lattice minimizes the Epstein zeta function among lattices of fixed density for a restricted range of real $s$. Cassels~\cite{Cassels1959} proposed an extension, but identified a gap in the proof in his corrigendum. Ennola~\cite{EnnolaLemma1964,EnnolaProblem1964} supplied a repair, and Diananda~\cite{Diananda1964} simplified the argument. Together, these works establish hexagonal-lattice minimality for the analytically continued zeta functions for real $s>0$, $s\ne1$. At $s=1$, the comparison concerns the finite difference between zeta functions of lattices of the same density, after cancellation of their common pole.

Montgomery~\cite{Montgomery1988} studied the lattice theta function
\begin{equation}
    \theta_\Lambda(\alpha) = \sum_{\lambda \in \Lambda} e^{-\pi \alpha \lVert \lambda \rVert^2}
    \qquad \text{for } \alpha > 0.
\end{equation}    
For $\operatorname{Re}s>1$, we have the Mellin transform identity, see, e.g., \cite{JorgensonLang2001}
\begin{equation}
    \zeta_\Lambda(s) = \frac{\pi^s}{\Gamma(s)} \int_{0}^{\infty} \left(\theta_\Lambda(\alpha) - 1\right) \alpha^{s-1} \, d\alpha,
\end{equation}
which relates the zeta function to the theta function. Montgomery proved that, among lattices of fixed density, the hexagonal lattice uniquely minimizes $\theta_\Lambda(\alpha)$ for every $\alpha>0$, up to orthogonal transformations. The Mellin identity therefore yields its optimality for the convergent Riesz energies. Integration against positive measures similarly gives minimality for completely monotone potentials of squared distance.

Beyond constrained optimization over rigid lattice structures, these geometric energy minimizers are deeply linked to the fundamental \emph{crystallization conjecture} in mathematical physics. This conjecture asserts that under appropriate physical conditions, a system of $N$ interacting particles will spontaneously break translation invariance and arrange itself into a periodic configuration as $N \to \infty$. As detailed in the comprehensive review by Blanc and Lewin~\cite{BlancLewin2015}, proving this phenomenon from first principles remains one of the most challenging open problems in statistical mechanics. Substantial progress in understanding the foundational limits of such systems has been achieved through the work of Serfaty~\cite{Serfaty2015} and Petrache and Serfaty~\cite{PetracheSerfaty2017}, who developed a rigorous variational framework for the next-order energy of Riesz gases, connecting macroscopic distribution profiles to microscopic point configurations. Bétermin~\cite{Betermin2016} established hexagonal-lattice optimality among planar Bravais lattices for completely monotone potentials of squared distance with finite energy at fixed density, and for certain Lennard--Jones-type potentials. Hardin and Saff~\cite{HardinSaff2005} studied energy asymptotics and limiting distributions of minimizing Riesz configurations on rectifiable manifolds.

Before the unconstrained problem was solved by OpenAI \cite{OpenAI}, optimality of the hexagonal lattice within certain classes of periodic configurations was obtained by Faulhuber, Shafkulovska, and Zlotnikov \cite{FauShaZlo24} and by Hardin and Tenpas \cite{HarTen25}. Moreover, recent work by Lebl\'e \cite{Leble2025} established that the hexagonal lattice is a local energy minimizer within periodic configurations.

\subsection{Disclosure on the usage of AI tools}
The authors used OpenAI's ChatGPT-6 Astra model and Microsoft's Copilot, based on the GPT-5.6 Sol (OpenAI) reasoning model, as interactive mathematical and computational assistants during the development of this work. The models were used to identify the necessary auxiliary results in \cite{OpenAI} needed to prove Theorem~\ref{thm:main-uniqueness}~(B). The formulation of the research questions, the assessment of the mathematical arguments, and the selection of results for inclusion remained with the authors. All claims in the final manuscript were reviewed, rewritten, and verified by the authors, who take full responsibility for the mathematical content and its presentation.

\section{Notation and further results}
Throughout, we denote the minimal distance $d$ between distinct points in $\Lambda_2$, defined by \eqref{eq:hexagonal-lattice}, by
\begin{equation}
    d=\sqrt{\frac{2}{\sqrt3}}.
\end{equation}

\begin{definition}
    For $n \in \N$, a periodic point configuration $C\subset\R^n$ is of the form
    \begin{equation}\label{eq:configuration}
        C=\bigcup_{j=1}^{N}(v_j+\Lambda),
        \qquad \rho(C)=\frac{N}{\vol(\R^n/\Lambda)},
    \end{equation}
    where $\Lambda\subset\R^n$ is a full-rank lattice and the $v_j$ are distinct modulo $\Lambda$, i.e., $v_i-v_j \notin \Lambda$ for $ i \neq j$.
\end{definition}
As $\vol(\R^n/\Lambda)$ is the volume of a fundamental cell of $\Lambda$, the periodic configuration $C$ in \eqref{eq:configuration} has density $1$ exactly when $\vol(\R^n/\Lambda)=N$.

For a nonnegative potential $p:\R_{>0} \to \R_{\geq 0}$ and a periodic configuration $C$ as in \eqref{eq:configuration}, the $p$-energy of $C$ is given by
\begin{equation}\label{eq:energy}
    E_p(C)=\frac1N\sum_{i,j=1}^{N}
    \sum_{\lambda\in\Lambda \setminus \{v_i-v_j\}}
    p\bigl(\lVert v_i-v_j+\lambda\rVert\bigr),
\end{equation}
where the self-interaction terms with $i=j$ and $\lambda=0$ are excluded. For $t>0$, put
\begin{equation}
    E_t(C)=E_{r \mapsto e^{-\pi t r^2}}(C).
\end{equation}
We call $E_t(C)$ the Gaussian energy of $C$. 
In the notation of the introduction, $E_t(C)=E_{\varphi_{\pi t}}(C)$. Thus \eqref{eq:energyminimality} applies with $\alpha=\pi t$.

We will prove Theorem~\ref{thm:main-uniqueness}~(A) in Section~\ref{sec:minimum-separation-route} and Theorem~\ref{thm:main-uniqueness}~(B) in Section~\ref{sec:auxiliary-function-route}. The arguments also yield the following additional results of independent interest. These rely on the connections between energy minimization and optimal sphere packings, which have been stated more informally, e.g., in \cite[Section~9]{CohnKumar2007}, \cite[Section~1.1]{CKMRV2022}, \cite[Section~1]{Radchenko2026}, and hold in arbitrary dimension.

\begin{proposition}\label{prop:uo-packing}
    For $n \in \N$, let $C\subset\R^n$ be a periodic configuration that minimizes $E_t$ among periodic configurations of the same density for every $t>0$. Then $C$ maximizes minimum separation among periodic configurations of that density. Equivalently, its centers give a densest periodic sphere packing. In particular, this conclusion holds for every universally optimal periodic configuration with density 1.
\end{proposition}

\begin{theorem}\label{thm:general-rigidity}
    For $n \in \N$, let $C\subset\R^n$ be a periodic configuration that minimizes $E_t$ among periodic configurations of the same density for every $t>0$. Suppose that, up to Euclidean motion, $C$ is the unique periodic configuration of that density maximizing minimum separation. Then, up to Euclidean motion, $C$ is the unique periodic configuration of that density minimizing every Gaussian energy simultaneously.
\end{theorem}

\section{Part I: the minimum separation route}\label{sec:minimum-separation-route}
The aim of this section is to prove Theorem~\ref{thm:main-uniqueness}~(A). 
We show that if a periodic configuration of density $1$ has the same energy as $\Lambda_2$ for every Gaussian potential, then it is congruent to $\Lambda_2$. 
The proof establishes the minimum separation of points in the configuration and uses the uniqueness of the densest periodic circle packing. Together with \eqref{eq:energyminimality}, this gives uniqueness of the simultaneous minimizer of all Gaussian energies.

We begin by showing that for periodic configurations, the decay rate of Gaussian energies for $t \to \infty$ corresponds directly to the minimum distance between points. For a periodic configuration $C$, the minimum distance between distinct points exists and is strictly positive. Indeed, a periodic configuration has only finitely many points modulo its period lattice. Its pair-displacement set is a finite union of translates of a discrete lattice and therefore contains only finitely many vectors in every bounded region. We denote this minimum distance by $\sepdist(C)$.

\begin{lemma}
    \label{lem:large-t}
    Let $C$ be a periodic configuration in $\R^n$. Then
    \begin{equation}\label{eq:log-limit}
        \lim_{t\to\infty} \left(-\frac1{\pi t}\log E_t(C) \right)=\sepdist(C)^2.
    \end{equation}
\end{lemma}

\begin{proof}
Write $r_0=\sepdist(C)$. Let $M>0$ be the average number of points at distance $r_0$ from a point of $C$, where the average is taken over the $N$ representatives $v_1,\ldots,v_N$ in one period cell, i.e.,
\begin{equation}
    M=\frac1N\sum_{i=1}^{N}|\{x\in C:\lVert x-v_i\rVert=r_0\}|.
\end{equation}

Then, by \eqref{eq:energy}
\begin{equation}\label{eq:lower-bound}
    E_t(C)\ge M e^{-\pi t r_0^2}.
\end{equation}
Fix $t_0>0$. For $t\ge t_0$, we use \eqref{eq:energy} and the fact that $\lVert v_i-v_j+\lambda\rVert\ge r_0$ in every summand to obtain
\begin{align}
    E_t(C)
    &=\frac1N\sum_{i,j=1}^{N}
    \sum_{\lambda\in\Lambda \setminus\{v_i-v_j\}}
    e^{-\pi t\lVert v_i-v_j+\lambda\rVert^2}\nonumber\\
    &=\frac{e^{-\pi t r_0^2}}{N}\sum_{i,j=1}^{N}
    \sum_{\lambda\in\Lambda \setminus\{v_i-v_j\}}
    e^{-\pi t(\lVert v_i-v_j+\lambda\rVert^2-r_0^2)}\nonumber\\
    &\le\frac{e^{-\pi t r_0^2}}{N}\sum_{i,j=1}^{N}
    \sum_{\lambda\in\Lambda \setminus\{v_i-v_j\}}
    e^{-\pi t_0(\lVert v_i-v_j+\lambda\rVert^2-r_0^2)} \nonumber\\
    &=e^{-\pi t r_0^2}e^{\pi t_0 r_0^2}E_{t_0}(C). \label{eq:upper-bound}
\end{align}

It is clear from \eqref{eq:energy} that the quantity $E_{t_0}(C)$ is finite for a periodic configuration. Taking logarithms in \eqref{eq:lower-bound} and \eqref{eq:upper-bound}, dividing by $-\pi t$, and letting $t\to\infty$ proves \eqref{eq:log-limit}.
\end{proof}

\begin{proof}[Proof of Proposition~\ref{prop:uo-packing}]
    Let $C$ satisfy the stated minimality assumption, and let $D$ be a periodic configuration of the same density. Hence,
    \begin{equation}
        E_t(C)\le E_t(D)
        \qquad \text{ for all } t>0.
    \end{equation}
    After taking logarithms, multiplying by $-1/(\pi t)$, and passing to the limit in Lemma~\ref{lem:large-t}, we obtain
    \begin{equation}
        \sepdist(C)^2\ge\sepdist(D)^2.
    \end{equation}
    Thus $C$ maximizes minimum separation. Balls of radius $\sepdist(C)/2$ centered at the elements of $C$ therefore form a densest periodic sphere packing. 
    A universally optimal periodic configuration of density $1$ satisfies the minimality assumption by definition.
\end{proof}

For periodic configurations of density $1$, Proposition~\ref{prop:uo-packing} gives the implication
\begin{equation}\label{eq:contrapositive}
    C \text{ is not an optimal sphere packing}
    \quad\Longrightarrow\quad
    C \text{ is not universally optimal}.
\end{equation}

\begin{proof}[Proof of Theorem~\ref{thm:general-rigidity}]
    Let $D$ be a periodic configuration of the same density as $C$, minimizing every Gaussian energy among periodic configurations of density $\rho(C)$. Since $C$ and $D$ both minimize $E_t$ for every $t>0$,
    \begin{equation}
        E_t(D)=E_t(C) \qquad \text{ for all } t>0.
    \end{equation}
    Lemma~\ref{lem:large-t} gives $\sepdist(D)=\sepdist(C)$. Thus, $D$ also maximizes minimum separation, and the assumed uniqueness of $C$ implies that $D$ is congruent to $C$.
\end{proof}

\begin{proof}[Proof of Theorem~\ref{thm:main-uniqueness}(A)]
    Let $C$ be a periodic configuration with $\rho(C) = 1$ for which equality holds in \eqref{eq:energyminimality} for every $\alpha>0$. Since $E_t(C)=E_{\varphi_{\pi t}}(C)$, this assumption gives
    \begin{equation}\label{eq:all-scales}
        E_t(C)=E_t(\Lambda_2),
        \qquad \text{ for all } t>0.
    \end{equation}
    Lemma~\ref{lem:large-t} yields
    \begin{equation}
        \sepdist(C)=\sepdist(\Lambda_2)=d.
    \end{equation}
    Hence disks of radius $d/2$ centered at the points of $C$ form a periodic packing. Since $\rho(C) = 1$, its packing density equals
    \begin{equation}\label{eq:packing-density}
        \pi\left(\frac d2\right)^2
        =\frac{\pi}{2\sqrt3}
        =\frac{\pi}{\sqrt{12}}.
    \end{equation}
    This is the maximal planar packing density. Since the hexagonal packing is the unique densest periodic circle packing in the plane up to Euclidean motion \cite{FejesToth1942, Thue1910} see also \cite{Boroczky2004, ConwaySloane1999, FejesToth1953}, it follows that $C$ is congruent to $\Lambda_2$.
\end{proof}

Two useful variants follow from the same argument.

\begin{corollary}
\label{cor:sequence}
    Let $C\subset\R^2$ be a periodic point configuration of density $1$, and let $\Lambda_2$ be the hexagonal lattice of density $1$ defined in \eqref{eq:hexagonal-lattice}. Let $(t_k)_{k\in\N}$ be a sequence of positive real numbers with $t_k\to\infty$. If
    \begin{equation}
        E_{t_k}(C)\le E_{t_k}(\Lambda_2)
    \end{equation}
    for every $k\in\N$, then $C$ is congruent to $\Lambda_2$.
\end{corollary}

\begin{proof}
Lemma~\ref{lem:large-t} gives $\sepdist(C)\ge d$.  Equation \eqref{eq:packing-density} and the uniqueness of the densest periodic circle packing in the plane completes the proof.
\end{proof}

\begin{corollary}
    \label{cor:unbounded-support}
    Let $g\not\equiv0$ be a completely monotone potential of squared distance with representation
    \begin{equation}\label{eq:bernstein}
        g(r)=\int_{(0,\infty)}e^{-\pi t r^2}\,d\mu(t)
    \end{equation}
    where $\mu$ is a positive Borel measure on $(0,\infty)$. Suppose $E_g(\Lambda_2)<\infty$ and $\operatorname{supp}(\mu)$ is unbounded. Then $\Lambda_2$ is the unique density-one periodic minimizer of $E_g$.
\end{corollary}

\begin{proof}
    Let $C\subset\R^2$ be a periodic configuration of density $1$. By Tonelli's theorem and \eqref{eq:energyminimality}, applied with $\alpha=\pi t$,
    \begin{equation}\label{eq:mixture}
        E_g(C)-E_g(\Lambda_2)
        =\int_{(0,\infty)}\bigl(E_t(C)-E_t(\Lambda_2)\bigr)\,d\mu(t)\ge0.
    \end{equation}
    If equality holds, then $E_t(C)=E_t(\Lambda_2)$ for $\mu$-almost every $t$. From the unbounded support assumption, we conclude that there is a diverging sequence $(t_k)_{k \in \N}$ with $E_{t_k}(C)=E_{t_k}(\Lambda_2)$, for all $k \in \N$. An application of Corollary~\ref{cor:sequence} completes the proof.
\end{proof}

\section{Part II: the auxiliary-function route}\label{sec:auxiliary-function-route}

We now use specific information from the proof in \cite{OpenAI}. This route yields the strong uniqueness statement of Theorem~\ref{thm:main-uniqueness}~(B). We use the Fourier transform
\begin{equation}
    \widehat f(\xi) = \int_{\R^2}f(x)e^{-2\pi i x\cdot\xi}\,dx, \quad \text{ for } \xi \in \R^2.
\end{equation}
It is well-known (see, e.g., \cite[Appendix~A]{Folland1989}) that the Fourier transform of a Gaussian $G_t$ is given by
\begin{equation}
    \widehat{G_t}(x) = \frac{1}{t} G_{1/t}(x),
    \quad
    G_t(x) = e^{-\pi t \lVert x \rVert^2},
    \qquad \text{ for all }
    x \in \R^2, \ t>0.
\end{equation}

We recall the following auxiliary-function construction from \cite[Lemma~6.2 and equation~(6.5)]{OpenAI}, expressed in our notation.

\begin{proposition}
    \label{prop:auxiliary}
    For every $t>0$ there is a real radial Schwartz function $f_t$ such that
    \begin{equation}\label{eq:auxiliary-signs}
        f_t \le G_t,
        \qquad
        \widehat{f}_t\ge0,
    \end{equation}
    and
    \begin{equation}\label{eq:sharpness}
        \widehat{f}_t(0)-f_t(0)=E_t(\Lambda_2).
    \end{equation}
    The functions $f_t$ and $\widehat f_t$ satisfy the contact conditions from \cite[Lemma~6.2]{OpenAI}:
    \begin{equation}
     f_t(x)=G_t(x)\quad \text{ for all } x\in\Lambda_2\setminus\{0\},
     \qquad
     \widehat{f}_t(\xi)=0\quad \text{ for all }\xi\in\Lambda_2^*\setminus\{0\},
    \end{equation}
    where $\Lambda_2^*$ denotes the dual lattice of $\Lambda_2$.
\end{proposition}
Recall that the dual lattice to $\Lambda$ is
\begin{equation}
    \Lambda^*=\{\xi\in\R^2:\xi\cdot\lambda\in\Z\text{ for every }\lambda\in\Lambda\}
\end{equation}

To prove Theorem~\ref{thm:main-uniqueness}~(B), we need $G_t(x)>f_t(x)$ whenever $0<\lVert x\rVert<d$. The following lemma establishes this property using the quantitative estimates in \cite{OpenAI}.

\begin{lemma}
\label{lem:strictness}
    For every $t>0$, choose $f_t$ as constructed in the proof of \cite[Lemma~6.2]{OpenAI}, and put 
    \begin{equation}
        Q_t=G_t-f_t
    \end{equation}
    Then
    
    \begin{equation}\label{eq:strictness}
        Q_t(x)>0
        \qquad \text{and} \qquad
        \widehat{f}_t(x) > 0,
        \qquad
        \text{ for all } 0<\lVert x\rVert<d.
    \end{equation}
\end{lemma}
\begin{proof}
    We first identify the functions in the construction of \cite{OpenAI}, then explain how its quantitative estimates imply strict positivity. We use the particular auxiliary functions specified in \cite{OpenAI} as the properties listed in Proposition~\ref{prop:auxiliary} alone do not give the quantitative estimates below.
    
    Throughout, we put
    \begin{equation}
        h=\frac25,\qquad
        k=\pi\left(d^2 t-h\right),\qquad
        s= \frac{\lVert x\rVert^2}{d^2}.
    \end{equation}
    
    \emph{Step 1: identify the functions for $t\ge1$.}
    As is common when working with Gaussian energies, we first treat the case $t\ge1$ and then employ the Fourier transform to cover the range $0<t<1$. The Gaussian parameter denoted by $\alpha$ in \cite[Lemma~6.2]{OpenAI} is our $t$. We have
    \begin{equation}
        k\ge\pi\left(\frac{2}{\sqrt3}-\frac25\right)>2.36,
    \end{equation}
    which is requested in the first sentence of \cite[Section 5]{OpenAI}.
    
    Also, we note the following equivalence
    \begin{equation}
        0<\lVert x\rVert<d\quad\Longleftrightarrow\quad 0<s<1.
    \end{equation}
    By \cite[Theorem 2.1]{OpenAI} it holds that there are entire functions $H_1$, $H_2$, real on the real axis, such that for all $x \in \R^2$
    \begin{equation}
        f(x) = e^{-\pi h \frac{\lVert x \rVert^2}{d^2}} H_1\left( \frac{\lVert x \rVert^2}{d^2} \right),
        \qquad
        \widehat{f}(x) 
        = e^{-\pi h \frac{\lVert x \rVert^2}{d^2}} H_2\left( \frac{\lVert x \rVert^2}{d^2} \right).
    \end{equation}
    As in \cite[Theorem~2.1]{OpenAI}, we set
    \begin{equation}
        T_k(s)=\frac{e^{-k(s-1)}}{k}.
    \end{equation}
    Next, we note that the proof of \cite[Lemma~6.2]{OpenAI} defines $f_t=ke^{-k}f$ in the third sentence.
    
    The relation $\pi h+k=\pi t d^2$ gives
    \begin{align}
        ke^{-k}e^{-\pi hs}T_k(s)
        &=e^{-k}e^{-\pi hs}e^{-k(s-1)}\nonumber\\
        &=e^{-(\pi h+k)s}
        =e^{-\pi t\lVert x\rVert^2}
        =G_t(x). \label{eq:alternativeFormOfGt}
    \end{align}
    The second sentence of Section~5 of \cite{OpenAI} defines
    \begin{equation}
        F_1=T_k-H_1,
        \qquad
        F_2=H_2.
    \end{equation}
    These notations imply that
    \begin{align}
        Q_t(x)
        &=G_t(x)-f_t(x)\nonumber\\
        &=ke^{-k}e^{-\pi hs}\bigl(T_k(s)-H_1(s)\bigr) \nonumber \\
        &=ke^{-k}e^{-\pi hs}F_1(s), \label{eq:auxiliary-normalizationQ}.
    \end{align}
    Linearity of the Fourier transform yields
    \begin{align}
        \widehat f_t(x)
        &=ke^{-k}\widehat f(x)
        =ke^{-k}e^{-\pi hs}F_2(s).\label{eq:auxiliary-normalizationK}
    \end{align}
    
    The common factor $ke^{-k}e^{-\pi hs}$ is strictly positive. 
    The assertion therefore follows if we prove $F_1(s),F_2(s)>0$ for $0<s<1$.
    
    \emph{Step 2: identify the two estimates needed.}
    We use the intervals $[0,1/2]$ and $[1/2,1]$. Equation~(5.1) of \cite{OpenAI}, applied at $1$, gives
    \begin{equation}\label{eq:normalized-contacts}
     F_i(1)=F_i'(1)=0,\qquad i=1,2.
    \end{equation}
    The quotient defined in \cite[Section~5.1]{OpenAI} therefore simplifies to
    \begin{equation}
     \mathcal D_{1,2}F_i(s)
     =\frac{F_i(s)-F_i(1)-(s-1)F_i'(1)}{(s-1)^2}
     =\frac{F_i(s)}{(s-1)^2}\qquad \text{ for } s\ne1.
    \end{equation}
    Therefore, to show positivity of $F_i$ for $0<s<1$, it suffices to prove positivity of $F_i$ on $[0,1/2]$ and of $\mathcal D_{1,2}F_i$ on $[1/2,1)$.
    
    To match the notation of \cite{OpenAI}, the two intervals are parametrized by
    \begin{align}\label{eq:mynu}
        s=m+yu,\qquad 0\le u\le1,\qquad
        (m,y,\nu)=\left(0,\frac12,0\right)
        \ \text{or}\
        \left(1,-\frac12,2\right).
    \end{align}
    The choice $\nu=0$ corresponds to estimating the function itself, while $\nu=2$ corresponds to estimating $\mathcal D_{1,2}F_i$. The second parametrization runs from $s=1$ to $s=1/2$.
    
    \emph{Step 3: apply the estimate from the proof of Lemma~5.2.}
    For
    \begin{equation}\label{eq:margin}
        \varepsilon
        =0.001-0.00011-0.000001-0.000001
        =0.000888>0,
    \end{equation}
    the displayed estimate near the end of the proof of \cite[Lemma~5.2]{OpenAI}, applied to the two choices of $(m,y,\nu)$ in \eqref{eq:mynu}, yields 
    \begin{equation}
        \begin{aligned}
            F_i(s)&>\varepsilon
            && \text{ for } 0\le s\le\frac12,\\
            \mathcal D_{1,2}F_i(s)&>\varepsilon
            && \text{ for } \frac12\le s<1,
        \end{aligned}
        \qquad i=1,2.
    \end{equation}
    We use this estimate from the proof because the statement of Lemma~5.2 asserts only nonnegativity.
    
    With Step 2 and Step 1, this yields \eqref{eq:strictness} for $t\ge1$.
    
    \emph{Step 4: treat $0<t<1$.}
    Put $\beta=1/t>1$. The functions $Q_\beta$ and $\widehat f_\beta$ are strictly positive for $0<\lVert x\rVert<d$ by the case already proved.
    Equation~(6.4) of \cite{OpenAI} defines the auxiliary function for $t$ by
    \begin{equation}
        f_t=\beta\,\widehat Q_\beta.
    \end{equation}
    With our Fourier convention, the Gaussian transform in dimension two is
    \begin{equation}
     \widehat{G_\beta}(x)
     =\frac1\beta e^{-\pi\lVert x\rVert^2/\beta}
     =\frac1\beta G_t(x).
    \end{equation}
    Since $Q_\beta=G_\beta-f_\beta$, it follows that
    \begin{equation}
        f_t=\beta(\widehat G_\beta-\widehat f_\beta)
        =G_t-\beta\widehat f_\beta,
        \qquad
        Q_t=G_t-f_t = \beta \widehat f_\beta.
    \end{equation}
    Moreover, applying the Fourier transform twice reflects a Schwartz function. The function $Q_\beta$ is radial and hence even. Therefore,
    \begin{equation}
     \widehat f_t(x)
           =\beta\,\widehat{\widehat {Q_\beta}}(x)
           =\beta Q_\beta(-x)
           =\beta Q_\beta(x).
    \end{equation}
    Both $Q_t=\beta\widehat f_\beta$ and $\widehat f_t=\beta Q_\beta$ are strictly positive on $0<\lVert x\rVert<d$, because $\beta>1>0$. Together with the case $t\ge1$, this proves the lemma.
\end{proof}

\begin{remark}\label{rem:certificate}
    Lemma~\ref{lem:strictness} makes no new computer-assisted estimate. It records the implication of the positive certified margin for the location of the first possible zero. The certificate, its error bounds, and Fourier complementation are all provided by \cite{OpenAI}.
\end{remark}

\begin{proof}[Proof of Theorem~\ref{thm:main-uniqueness}~(B)]
    Let $C\subset\R^2$ be a periodic configuration with $\rho(C) = 1$ for which equality holds in \eqref{eq:energyminimality} at some $\alpha_0>0$. Set $t=\alpha_0/\pi$, so that $E_t(C)=E_t(\Lambda_2)$. Choose the auxiliary function $f_t$ from Proposition~\ref{prop:auxiliary}, with the strictness property established in Lemma~\ref{lem:strictness}. Write $C$ as in \eqref{eq:configuration};
    \begin{equation}
        C = \bigcup_{j=1}^N (v_j+\Lambda),
        \qquad \Lambda \text{ a full-rank lattice in } \R^2.
    \end{equation}
    Then $\rho(C) = 1$ implies $\vol(\R^2/\Lambda)=N$. To keep the proof self-contained, we repeat the arguments leading to the linear programming bounds for energy minimization \cite{CohnKumar2007}, which, just as the bounds for sphere packings \cite{CohnElkies2003}, follow from the Poisson summation formula. Define
    \begin{equation}
        S(\xi)=\sum_{j=1}^{N}e^{2\pi i\xi\cdot v_j}.
    \end{equation}
    Since $f_t$ is a Schwartz function, the Poisson summation formula applies with absolute convergence of both sums. With our Fourier-transform convention and $\vol(\R^2/\Lambda)=N$, it yields
    \begin{equation}
     \sum_{\lambda\in\Lambda}f_t(v_i-v_j+\lambda)
     =\frac1N\sum_{\xi\in\Lambda^*}
     \widehat{f}_t(\xi)e^{2\pi i\xi\cdot(v_i-v_j)}.
    \end{equation}
    Summing over $i,j$ and using
    \begin{equation}
     \sum_{i,j=1}^{N}e^{2\pi i\xi\cdot(v_i-v_j)}=|S(\xi)|^2
    \end{equation}
    yields
    \begin{equation}
     \frac1N\sum_{i,j=1}^{N}\sum_{\lambda\in\Lambda}
     f_t(v_i-v_j+\lambda)
     =\frac1{N^2}\sum_{\xi\in\Lambda^*}\widehat{f}_t(\xi)|S(\xi)|^2.
    \end{equation}
    The energy in \eqref{eq:energy} omits the $N$ diagonal terms with $i=j$ and $\lambda=0$, which contribute $f_t(0)$ after normalization. Also note that $S(0)=N$. Hence, by \eqref{eq:sharpness} 
    \begin{align}
        \frac1N\sum_{i,j=1}^{N}
        \sum_{\lambda\in\Lambda \setminus \{v_i-v_j\}}
        f_t(v_i-v_j+\lambda) &= \frac1{N^2}\sum_{\xi\in\Lambda^*\setminus\{0\}}
        \widehat{f}_t(\xi)|S(\xi)|^2 + \widehat{f}_t(0) - f_t(0) \nonumber\\
        &= \frac1{N^2}\sum_{\xi\in\Lambda^*\setminus\{0\}}
        \widehat{f}_t(\xi)|S(\xi)|^2 + E_t(\Lambda_2).\label{eq:sharpnessapplied}
    \end{align}
    Since $G_t=f_t+Q_t$ we conclude from \eqref{eq:sharpnessapplied}
    \begin{align}
        E_t(C)-E_t(\Lambda_2)
        &=\frac1N\sum_{i,j=1}^{N}
        \sum_{\lambda\in\Lambda \setminus \{v_i-v_j\}}
        Q_t(v_i-v_j+\lambda) + \frac1N\sum_{i,j=1}^{N}
        \sum_{\lambda\in\Lambda \setminus \{v_i-v_j\}}
        f_t(v_i-v_j+\lambda) -E_t(\Lambda_2)\nonumber\\
        &= \frac1N\sum_{i,j=1}^{N}
        \sum_{\lambda\in\Lambda \setminus \{v_i-v_j\}}
        Q_t(v_i-v_j+\lambda) +\frac1{N^2}\sum_{\xi\in\Lambda^*\setminus\{0\}}
        \widehat{f}_t(\xi)|S(\xi)|^2.\label{eq:lp-identity}
    \end{align}
    Every summand in \eqref{eq:lp-identity} is nonnegative by \eqref{eq:auxiliary-signs}. If $E_t(C)=E_t(\Lambda_2)$, both sums on the right-hand side of \eqref{eq:lp-identity} are zero. In particular, each nonnegative term in the first sum must vanish:
    \begin{equation}
        Q_t(v_i-v_j+\lambda)=0
    \end{equation}
    for every $i,j\in\{1,\ldots,N\}$ and $\lambda\in\Lambda$ with $i\ne j$ or $\lambda\ne0$. By \eqref{eq:strictness}, every such nonzero vector $v_i-v_j+\lambda$ has length at least $d$.
    Thus $\sepdist(C) \ge d$, and the packing argument used in the proof of Theorem~\ref{thm:main-uniqueness}~(A) gives $C \cong \Lambda_2$.
\end{proof}

\begin{proof}[Proof of Corollary~\ref{cor:cm}]
    By definition, we have $g(r)=f(r^2)$ for some nonzero completely monotone function $f$. By Theorem~\ref{thm:bernstein-widder},
    \begin{equation}
        f(s)=\int_{[0,\infty)}e^{-\pi ts}\,d\mu(t)
    \end{equation}
    for a nonzero positive measure $\mu$. If $\mu(\{0\})>0$, then $g(r)=f(r^2)\ge\mu(\{0\})>0$ for every $r>0$, contradicting the finite-energy assumption. Thus $\mu(\{0\})=0$, and
    \begin{equation}
        g(r)=f(r^2)=\int_{(0,\infty)}e^{-\pi tr^2}\,d\mu(t).
    \end{equation}
    
    Let $C\subset\R^2$ be periodic of density $1$. Tonelli's theorem and \eqref{eq:energyminimality}, applied with $\alpha=\pi t$, give
    \begin{equation}
        E_g(C)-E_g(\Lambda_2)
        =\int_{(0,\infty)}\bigl(E_t(C)-E_t(\Lambda_2)\bigr)\,d\mu(t)\ge0.
    \end{equation}
    If $E_g(C)=E_g(\Lambda_2)$, then $E_t(C)=E_t(\Lambda_2)$ for $\mu$-almost every $t$. Since $\mu\ne0$, equality holds for at least one $t>0$, and Theorem~\ref{thm:main-uniqueness}(B), applied with $\alpha=\pi t$, gives $C \cong \Lambda_2$.
\end{proof}


\begin{thebibliography}{99}

\bibitem{Bernstein1929}
S.~Bernstein,
Sur les fonctions absolument monotones,
\emph{Acta Math.} 52:1--66, 1929.
\doi{10.1007/BF02592679}.

\bibitem{Betermin2016}
L.~Bétermin,
Two-dimensional theta functions and crystallization among Bravais lattices,
\emph{SIAM J. Math. Anal.} 48(5):3236--3269, 2016.
\doi{10.1137/15M101614X}.

\bibitem{BlancLewin2015}
X.~Blanc and M.~Lewin, The crystallization conjecture: a review,
\emph{EMS Surv. Math. Sci.} 2(2):255--306, 2015.
\doi{10.4171/EMSS/13}.

\bibitem{Boroczky2004}
K.~B\"or\"oczky Jr.,
\emph{Finite Packing and Covering},
Cambridge Tracts in Mathematics, vol.~154,
Cambridge University Press, Cambridge, 2004.
\doi{10.1017/CBO9780511546587}.

\bibitem{Cassels1959}
J.~W.~S. Cassels,
On a problem of Rankin about the Epstein zeta-function,
\emph{Glasg. Math. J.} 4(2):73--80, 1959.
\doi{10.1017/S2040618500033906}.
(Corrigendum, ibid. 6:116, 1963.
\doi{10.1017/S2040618500034833}.)

\bibitem{CohnElkies2003}
H.~Cohn and N.~Elkies,
New upper bounds on sphere packings I,
\emph{Ann. of Math.} (2) \textbf{157} (2003), no.~2, 689--714.
\doi{10.4007/annals.2003.157.689}.

\bibitem{CohnKumar2007}
H.~Cohn and A.~Kumar,
Universally optimal distribution of points on spheres,
\emph{J. Amer. Math. Soc.} 20(1):99--148, 2007.
\doi{10.1090/S0894-0347-06-00546-7}.

\bibitem{CKMRV2017}
H.~Cohn, A.~Kumar, S.~D. Miller, D.~Radchenko, and M.~Viazovska,
The sphere packing problem in dimension \(24\),
\emph{Ann. of Math. (2)} 185(3):1017--1033, 2017.
\doi{10.4007/annals.2017.185.3.8}.

\bibitem{CKMRV2022}
H.~Cohn, A.~Kumar, S.~D. Miller, D.~Radchenko, and M.~Viazovska,
Universal optimality of the \(E_8\) and Leech lattices and interpolation formulas,
\emph{Ann. of Math. (2)} 196(3):983--1082, 2022.
\doi{10.4007/annals.2022.196.3.3}.

\bibitem{ConwaySloane1999}
J.~H. Conway and N.~J.~A. Sloane,
\emph{Sphere Packings, Lattices and Groups},
3rd ed., Grundlehren der Mathematischen Wissenschaften, vol.~290,
Springer-Verlag, New York, 1999.
\doi{10.1007/978-1-4757-6568-7}.

\bibitem{Diananda1964}
P.~H. Diananda,
Notes on two lemmas concerning the Epstein zeta-function,
\emph{Glasg. Math. J.} 6(4):202--204, 1964.
\doi{10.1017/S2040618500035036}.

\bibitem{EnnolaLemma1964}
V.~Ennola,
A lemma about the Epstein zeta-function,
\emph{Glasg. Math. J.} 6(4):198--201, 1964.\\
\doi{10.1017/S2040618500035024}.

\bibitem{EnnolaProblem1964}
V.~Ennola,
On a problem about the Epstein zeta-function,
\emph{Math. Proc. Cambridge Philos. Soc.} 60(4):855--875, 1964.
\doi{10.1017/S0305004100038330}.

\bibitem{FauShaZlo24}
M. Faulhuber, I. Shafkulovska, and I. Zlotnikov,
A note on energy minimization in dimension 2,
\emph{Proc. Amer. Math. Soc. Ser. B} 11:664--679, 2024.
\doi{10.1090/bproc/247}.

\bibitem{FejesToth1942}
L.~Fejes T\'oth,
Über die dichteste Kugellagerung,
\emph{Math. Z.}, 48:676--684, 1942.
\doi{10.1007/BF01180035}.

\bibitem{FejesToth1953}
L.~Fejes T\'oth,
\emph{Lagerungen in der Ebene, auf der Kugel und im Raum},
Die Grundlehren der Mathematischen Wissenschaften, vol.~65,
Springer-Verlag, Berlin, 1953.
\doi{10.1007/978-3-662-01206-2}.

\bibitem{Folland1989}
G.~B.~Folland,
\emph{Harmonic analysis in phase space},
Annals of Mathematics Studies, number~122,
Princeton University Press, 1989.
\doi{10.1515/9781400882427}.

\bibitem{HardinSaff2005}
D.~P. Hardin and E.~B. Saff,
Minimal Riesz energy point configurations for rectifiable $d$-dimensional manifolds,
\emph{Adv. Math.} 193(1):174--204, 2005.
\doi{10.1016/j.aim.2004.05.006}.

\bibitem{HarTen25}
D.~P. Hardin and N. Tenpas,
Universally Optimal Periodic Configurations in the Plane,
\emph{Discrete Anal.} 26, 2025.
\doi{10.48550/arXiv.2307.15822}.

\bibitem{JorgensonLang2001}
J.~Jorgenson und S.~Lang,
\emph{The Ubiquitous Heat Kernel},
in: B.~Engquist und W.~Schmid (Hrsg.),
\emph{Mathematics Unlimited --- 2001 and Beyond},
Springer, Berlin, Heidelberg, 2001, S.~655--683.
\doi{10.1007/978-3-642-56478-9\_34}.

\bibitem{Leble2025}
T.~Lebl{\'e},
The hexagonal lattice is universally locally optimal,
\emph{arXiv preprint} arXiv:2511.03353, 2025.
\doi{10.48550/arXiv.2511.03353}.

\bibitem{Montgomery1988}
H.~L. Montgomery,
Minimal theta functions,
\emph{Glasgow Math. J.} 30(1):75--85, 1988.
\doi{10.1017/S0017089500007047}.

\bibitem{OpenAI}
OpenAI,
\emph{Universal optimality of the triangular lattice},
preprint, 2026.
Available at:\\
\url{https://github.com/openai/math/blob/main/preprints/Universal-optimality-of-the-triangular-lattice-September-23-2026}

\bibitem{PetracheSerfaty2017}
M.~Petrache and S.~Serfaty,
Next order asymptotics and renormalized energy for Riesz interactions,
\emph{J. Inst. Math. Jussieu} 16(3):501--569, 2017.
\doi{10.1017/S1474748015000201}.

\bibitem{Radchenko2026}
D.~Radchenko,
\emph{Fourier uniqueness pairs and interpolation formulas},
in: G.~P.~Curbera, M.~L.~Puertas, A.~Su{\'a}rez-Llorens (Eds.),
European Congress of Mathematics (Sevilla, July 15--19, 2024),
EMS Press, pp. 975--990, 2026.
\doi{10.4171/9ECM/40}.

\bibitem{Rankin1953}
R.~A. Rankin,
A minimum problem for the Epstein zeta-function,
\emph{Glasg. Math. J.} 1(4):149--158, 1953.
\doi{10.1017/S2040618500035668}.

\bibitem{Serfaty2015}
S.~Serfaty,
\emph{Coulomb Gases and Ginzburg-Landau Vortices},
Zurich Lectures in Advanced Mathematics, European Mathematical Society (EMS), Zürich, 2015.
\doi{10.4171/152}.

\bibitem{Thue1910}
A. Thue,
Über die dichteste Zusammenstellung von kongruenten Kreisen in einer Ebene.
\emph{Norske Vid. Selsk. Skr.} 1:1--9, 1910.
Available at:
\url{https://www.biodiversitylibrary.org/page/15174094}

\bibitem{VentevogelNijboer1979}
W.~J. Ventevogel and B.~R.~A. Nijboer,
On the configuration of systems of interacting particle with minimum potential energy per particle,
\emph{Physica A} 98(1):274--288, 1979.
\doi{10.1016/0378-4371(79)90178-X}.

\bibitem{Viazovska2017}
M.~S. Viazovska,
The sphere packing problem in dimension \(8\),
\emph{Ann. of Math. (2)} 185(3):991--1015, 2017.
\doi{10.4007/annals.2017.185.3.7}.

\bibitem{Widder1941}
D.~V. Widder,
\emph{The Laplace Transform},
Princeton Mathematical Series, vol.~6,
Princeton University Press, Princeton, 1941.
\doi{10.1515/9781400876457}.

\end{thebibliography}
\end{document}